\documentclass[10pt,a4paper]{article}
\usepackage{amsmath, amssymb}
\usepackage{latexsym, color}
\usepackage{epsfig}

\numberwithin{equation}{section}

\newtheorem{thm}{Theorem}[section]

\newtheorem{lem}[thm]{Lemma}
\newtheorem{prop}[thm]{Proposition}

\def\nm{\noalign{\medskip}}

\newcommand{\qed}{\hfill \ensuremath{\square}}
\newcommand{\ds}{\displaystyle}
\newcommand{\pf}{\noindent {\sl Proof}. \ }
\newcommand{\p}{\partial}
\newcommand{\pd}[2]{\frac {\p #1}{\p #2}}

\newcommand{\eqnref}[1]{(\ref {#1})}

\newcommand{\Ibb}{\mathbb{I}}

\newcommand{\Kbb}{\mathbb{K}}
\newcommand{\Rbb}{\mathbb{R}}

\newcommand{\la}{\langle}
\newcommand{\ra}{\rangle}

\newcommand{\Hcal}{\mathcal{H}}

\newcommand{\Kcal}{\mathcal{K}}

\newcommand{\Scal}{\mathcal{S}}

\newcommand{\Gd}{\delta}
\newcommand{\Ge}{\epsilon}

\newcommand{\Gvf}{\varphi}

\newcommand{\Gl}{\lambda}

\newcommand{\GD}{\Delta}

\newcommand{\GG}{\Gamma}

\newcommand{\GO}{\Omega}

\newcommand{\beq}{\begin{equation}}
\newcommand{\eeq}{\end{equation}}

\begin{document}
\title{Spectral theory of a Neumann-Poincar\'e-type operator and analysis of
anomalous localized resonance II\thanks{\footnotesize This work was
supported by the ERC Advanced Grant Project MULTIMOD--267184 and
NRF grants No. 2010-0004091, and 2010-0017532, and by the NSF through grants DMS-0707978 and DMS-1211359.}}

\author{Habib Ammari\thanks{\footnotesize Department of Mathematics and Applications, Ecole Normale Sup\'erieure,
45 Rue d'Ulm, 75005 Paris, France (habib.ammari@ens.fr).} \and
Giulio Ciraolo\thanks{\footnotesize  Dipartimento di Matematica e
Informatica, Universit\`a di Palermo Via Archirafi 34, 90123, Palermo, Italy
  (g.ciraolo@math.unipa.it).} \and Hyeonbae
Kang\thanks{Department of Mathematics, Inha University, Incheon
402-751, Korea (hbkang@inha.ac.kr, hdlee@inha.ac.kr).}  \and
Hyundae Lee\footnotemark[4]  \and Graeme W.
Milton\thanks{\footnotesize Department of Mathematics, University
of Utah, Salt Lake City, UT 84112, USA (milton@math.utah.edu).}}

\maketitle

\begin{abstract}
If a body of dielectric material is coated by a plasmonic
structure of negative dielectric constant with nonzero loss
parameter, then cloaking by anomalous localized resonance (CALR)
may occur as the loss parameter tends to zero. The aim of this
paper is to investigate this phenomenon in two and three
dimensions when the coated structure is radial, and the core, shell and
matrix are isotropic materials. In two dimensions,
we show that if the real part of the permittivity of the shell is
$-1$ (under the assumption that the permittivity of the background
is $1$), then CALR takes place. If it is different from $-1$, then
CALR does not occur. In three dimensions, we show that CALR does
not occur. The analysis of this paper reveals that occurrence of
CALR is determined by the eigenvalue distribution of the
Neumann-Poincar\'e-type operator associated with the structure.
\end{abstract}

%%%%%%%%%%%%%%%%%%%%%%%%%%%%%
\section{Introduction}
%%%%%%%%%%%%%%%%%%%%%%%%%%%%%

If a body of dielectric material is coated by a plasmonic structure of negative dielectric constant
(with nonzero loss parameter), then anomalous localized resonance may occur as the loss parameter
tends to zero. This phenomena, first discovered by Nicorovici, McPhedran and Milton \cite{NMM_94} (see also \cite{MNMP_PRSA_05}), is responsible for the
subwavelength focussing properties of superlenses \cite{superlens}, and also occurs in magnetoelectric and
thermoelectric systems \cite{MNMP_PRSA_05}. The fields blow-up in a localized region, which moves as the position
of the source is moved, which is why it is termed anomalous localized resonance.
Remarkably, as found by Milton and Nicorovici \cite{MN_PRSA_06}
the localized resonant fields created by a source can act back on the source and cloak it. This invisibility cloaking
has attracted much attention \cite{MN_PRSA_06,osa,bruno,MNM_07,NMET_08,MNMCJ_09,MNBM_09,bouchitte, Dong, NMB_11, acklm, klsw, xiao}.

To state the problem and results in a precise way, let $\GO$ be a bounded domain in $\Rbb^d$, $d=2,3$, and $D$
be a domain whose closure is contained in $\GO$. For a given loss parameter $\Gd>0$, the permittivity
distribution in $\Rbb^d$ is given by
 \beq
 \Ge_\Gd  = \begin{cases}
 1 \quad & \mbox{in } \Rbb^d \setminus \overline{\GO}, \\
\Ge_s+ i \Gd \quad & \mbox{in } \GO \setminus \overline{D}, \\
 \Ge_c \quad &\mbox{in } D,
 \end{cases}
 \eeq
where $-\Ge_s$ and $\Ge_c$ are positive. We may consider the configuration as a core with permittivity $\Ge_c$ coated by the shell $\GO \setminus \overline{D}$ with permittivity $\Ge_s+i\Gd$.
For a given function $f$ compactly supported in $\Rbb^d \setminus \overline{\GO}$ satisfying
 \beq\label{zeroint}
 \int_{\Rbb^2} f\, dx=0
 \eeq
(which is required by conservation of charge), we consider the following dielectric problem:
 \beq \label{basiceqn}
 \nabla \cdot \Ge_\Gd \nabla V_\Gd = f \quad \mbox{in } \Rbb^d,
 \eeq
with the decay condition $V_\Gd (x) \to 0$ as $|x| \to \infty$.
The problem of cloaking by anomalous localized resonance (CALR) can be
formulated as the problem of identifying the sources $f$ such that
 first
 \beq\label{blowup1}
E_\Gd := \int_{\GO\setminus D} \Gd |\nabla V_\Gd|^2\, dx \to \infty \quad\mbox{as } \Gd \to 0,
 \eeq
and second, $V_\Gd/\sqrt{E_\Gd}$ goes to zero outside some radius $a$, as $\Gd\to 0$:
\beq\label{bounded}
|V_\Gd (x)/\sqrt{E_\Gd}| \to 0  \quad\mbox{as } \Gd \to 0 \quad\mbox{when}\, |x| > a.
\eeq
Physically the quantity $E_\Gd$ is proportional to the electromagnetic power dissipated into heat by the time harmonic electrical field averaged
over time. Using integration by parts we have the identity
\beq\label{conserv}
E_\Gd=\Im \int_{\Rbb^d} (\Ge_\Gd \nabla V_\Gd)\cdot\nabla \overline{V_\Gd}\, dx
=- \Im \int_{\Rbb^d} f\overline{V_\Gd}\, dx
\eeq
which equates the power dissipated into heat with the electromagnetic power produced by the source, where $\overline{V_\Gd}$ is the complex conjugate of $V_\Gd$.
Hence \eqnref{blowup1} implies an infinite amount of energy dissipated per unit time in the
limit $\Gd\to 0$ which is unphysical. If we rescale the source $f$ by a factor of $1/\sqrt{E_\Gd}$
then the source will produce the same power independent of $\Gd$ and the new associated potential $V_\Gd/\sqrt{E_\Gd}$ will,
by  \eqnref{bounded}, approach zero outside the radius $a$: cloaking due
to anomalous localized resonance (CALR) occurs.

In the recent paper \cite{acklm} the authors develop a spectral
approach to analyze the CALR phenomenon. In particular, they
show that if $D$ and $\GO$ are concentric disks in $\Rbb^2$ and
$\Ge_c=-\Ge_s=1$, then there is a critical radius $r_*$ such that
for any source $f$ supported outside $r_*$ CALR does not occur,
and for sources $f$ satisfying a mild condition CALR takes place.
The critical radius $r_*$ is given by
\begin{equation} \label{rcrid2}
r_* = \sqrt{r_e^3/r_i}, \end{equation}  where $r_e$ and $r_i$ are
the radii of $\GO$ and $D$, respectively. It is worth mentioning
that these results were extended in \cite{klsw} to the case when
the core $D$ is not radial by a different method based on a
variational approach.

The purpose of this paper is to extend some of the results in
\cite{acklm} in two directions. We consider the case when $\Ge_c$
and $-\Ge_s$ are not both $1$ and we consider CALR in three dimensions. The results
of this paper are threefold: Let $\GO$ and $D$ be concentric disks
or balls in $\Rbb^d$ of radii $r_e$ and $r_i$, respectively. Then,
the following results hold:
\begin{itemize}
\item If $d=2$ and $\Ge_s=-1$, then CALR occurs. When $\Ge_c=1$ the critical radius  $r_*$ is given
by \eqnref{rcrid2} and when $\Ge_c \neq 1$ the critical radius is
\beq r_*= \frac{r_e^2}{r_i}. \eeq
That is, for almost any source $f$
supported inside $r_*$ CALR occurs and for any source $f$
supported outside $r_*$ CALR does not occur. When $\Ge_c \neq 1$ the cloaking radius $r_e^2/r_i$ matches
that found in \cite{MN_PRSA_06}  for a single dipolar source (see figure 5 there and accompanying text).
\item If $\Ge_s \neq -1$, then CALR does not occur.
\item If $d=3$, then CALR does not
occur whatever $\Ge_s$ and $\Ge_c$ are.
\end{itemize}

We emphasize that the result on non-occurrence of CALR in three dimensions holds only when the dielectric constant $\Ge_s$ is constant. In the recent work \cite{acklm3} we show that CALR does occur in three dimensions if we use a shell with non-constant (anisotropic) dielectric constant.

It turns out that the occurrence of CALR depends on the
distribution of eigenvalues of the Neumann-Poincar\'e (NP)
operator associated with the structure (see the next section for
the definition of the NP operator). The NP operator is compact with its
eigenvalues accumulating towards $0$. It is proved in \cite{acklm} that
in two dimensions the NP operator associated with the circular
structure has the eigenvalues $\pm \rho^n$ for $n=1,2, \ldots$,
where $\rho=r_i/r_e$. We show that in three dimensions the NP
operator associated with the spherical structure has the
eigenvalues
\beq
\pm \frac{1}{2(2n+1)}\sqrt{1+4n(n+1)\rho^{2n+1}},
\quad n=0,1,\ldots.
\eeq
The exponential convergence of the
eigenvalues in two dimensions is responsible for the occurrence of
CALR and the slow convergence (at the rate $1/n$) in three dimensions
is responsible for the non-occurrence.

%%%%%%%%%%%%%%%%%%%%%%%%%%%%%%%%%%%%%%%%%%%%%
%
\section{Layer potential formulation}
%
%%%%%%%%%%%%%%%%%%%%%%%%%%%%%%%%%%%%%%%%%%%%%

Let $G$ be the fundamental solution to the Laplacian in $\Rbb^d$ which is given by
$$
G(x) =
\left\{
\begin{array}{ll}
\ds \frac{1}{2\pi} \ln |x|, \quad & d=2 , \\
\nm
\ds -\frac{1}{4\pi} \frac{1}{|x|}, \quad & d=3.
\end{array}
\right.
$$
Let $\Gamma_i:= \p D$ and $\Gamma_e:= \p\GO$. For
$\Gamma=\Gamma_i$ or $\Gamma_e$, we denote the single layer
potential of a function $\Gvf \in L^2(\Gamma)$ as
$\Scal_\Gamma[\Gvf]$, where
\begin{align*}
\Scal_\Gamma[\Gvf] (x) &:= \int_{\Gamma} G(x-y) \Gvf (y) \, d
\sigma(y), \quad x \in \Rbb^d.
\end{align*}
We also define the boundary integral operator $\Kcal_\Gamma$ on $L^2(\Gamma)$ by
 $$
 \Kcal_\Gamma[\Gvf] (x) := \int_{\Gamma} \frac{\p G(x-y)}{\p\nu(y)} \Gvf (y)\,d\sigma(y), \quad x \in \GG, 
 $$
and let
$\Kcal_\Gamma^*$ be the  $L^2$-adjoint of $\Kcal_\Gamma$. Hence,
the operator $\Kcal_\Gamma^*$ is given by
$$
\Kcal_\Gamma^*[\Gvf](x) = \int_{\Gamma} \frac{\p G(x-y)}{\p\nu(x)} \Gvf (y)\,d\sigma(y), \quad \Gvf \in  L^2(\Gamma).
 $$
The operators $\Kcal_\Gamma$ and
$\Kcal_\Gamma^*$ are sometimes called Neumann-Poincar\'e
operators. These operators are compact in $L^2(\Gamma)$ if
$\Gamma$ is $\mathcal{C}^{1, \alpha}$ for some $\alpha>0$.

The following notation will be used throughout this paper. For a
function $u$ defined on $\Rbb^d \setminus \Gamma$, we denote
 $$
 u |_\pm(x) := \lim_{t \to 0^+}  u(x \pm t \nu(x)), \quad x \in \Gamma ,
 $$
and
 $$
 \pd{u}{\nu} \Big |_\pm(x) := \lim_{t \to 0^+} \la \nabla u(x \pm t \nu(x)), \nu(x) \ra\;, \quad x \in \Gamma,
 $$
if the limits exist. Here and throughout this paper, $\la \;, \;\ra$ denotes the scalar
product on $\Rbb^d$.

The following jump formula relates the traces of the normal derivative of the single layer potential
to the operator $\Kcal_\Gamma^*$. We have
 \begin{align}
 \pd{}{\nu} \Scal_\Gamma [\Gvf] \Big |_\pm (x) & = \left( \pm
 \frac{1}{2} I + \Kcal_\Gamma^* \right) [\Gvf] (x), \quad
 x \in \Gamma. \label{singlejump}
 \end{align}
Here, $\nu$ is the outward unit normal vector field to $\Gamma$. See, for example, \cite{book2, Folland76}.

Let $F$ be the Newtonian potential of $f$, {\it i.e.},
 \beq \label{newton}
 F(x)= \int_{\Rbb^d} G(x-y) f(y) dy, \quad x \in \Rbb^d.
 \eeq
Then $F$ satisfies $\Delta F=f$ in $\Rbb^d$, and the solution $V_\Gd$ to \eqnref{basiceqn} may be
represented as
 \beq\label{vdelta}
 V_\Gd(x) = F(x) + \Scal_{\Gamma_i}[\Gvf_i](x) + \Scal_{\Gamma_e}[\Gvf_e](x)
 \eeq
for some functions $\Gvf_i \in L^2_0(\Gamma_i)$ and $\Gvf_e \in
L^2_0(\Gamma_e)$ ($L^2_0$ is the collection of all square
integrable functions with zero mean-value). The transmission
conditions along the interfaces $\Gamma_e$ and $\Gamma_i$
satisfied by $V_\Gd$ read
 \begin{align*}
 (\Ge_s+i\Gd) \pd{V_\Gd}{\nu} \Big|_{+} =  \Ge_c\pd{V_\Gd}{\nu}
 \Big|_{-} \quad \mbox{on } \Gamma_i , \\
 \pd{V_\Gd}{\nu} \Big|_{+} = (\Ge_s+i\Gd) \pd{V_\Gd}{\nu} \Big|_{-} \quad \mbox{on }
 \Gamma_e.
 \end{align*}
Hence the pair of potentials $(\Gvf_i, \Gvf_e)$  is the solution to
the following system of integral equations:
 $$
 \begin{cases}
 \ds (\Ge_s+i\Gd) \pd{\Scal_{\Gamma_i}[\Gvf_i]}{\nu_i} \Big|_{+} -
  \Ge_c\pd{\Scal_{\Gamma_i}[\Gvf_i]}{\nu_i} \Big|_{-} +
 (\Ge_s-\Ge_c+i\Gd) \pd{\Scal_{\Gamma_e}[\Gvf_e]}{\nu_i}  = (-\Ge_s+\Ge_c-i\Gd) \pd{F}{\nu_i} \quad \mbox{on } \Gamma_i, \\
 \nm
 \ds (-1+\Ge_s +i\Gd) \pd{\Scal_{\Gamma_i}[\Gvf_i]}{\nu_e} -
 \pd{\Scal_{\Gamma_e}[\Gvf_e]}{\nu_e} \Big|_{+} +
 (\Ge_s+i\Gd) \pd{\Scal_{\Gamma_e}[\Gvf_e]}{\nu_e} \Big|_{-}
 = (1-\Ge_s-i\Gd) \pd{F}{\nu_e} \quad \mbox{on } \Gamma_e.
 \end{cases}
 $$
Note that we have used the notation $\nu_i$ and $\nu_e$ to indicate
the outward normal on $\Gamma_i$ and $\Gamma_e$, respectively.
Using the jump formula \eqnref{singlejump} for the normal derivative of the single
layer potentials, the above equations can be rewritten as
\begin{equation}\label{matrixeq3}
\begin{bmatrix}
 \ds z_i^\Gd  I - \Kcal_{\Gamma_i}^* & & \ds - \pd{}{\nu_i} \Scal_{\Gamma_e} \\
 \nm \ds \pd{}{\nu_e} \Scal_{\Gamma_i} & & z_e^\Gd  I + \ds \Kcal_{\Gamma_e}^*
 \end{bmatrix}  \begin{bmatrix} \ds \Gvf_i \\ \Gvf_e \end{bmatrix} = \begin{bmatrix} \ds
 \frac{\p F}{\p \nu_i} \\ \nm \ds - \frac{\p F}{\p \nu_e} \end{bmatrix}
\end{equation}
on $\Hcal_0=L^2_0(\Gamma_i) \times L^2_0(\Gamma_e)$, where we set
\begin{equation}\label{lam e mu}
z_i^\Gd  = \frac{\Ge_c+\Ge_s+i\Gd}{2(\Ge_c-\Ge_s-i\Gd)}, \quad z_e^\Gd  = \frac{1+\Ge_s+i\Gd}{2(1-\Ge_s-i\Gd)}.
\end{equation}

Let $\Hcal=L^2(\Gamma_i) \times L^2(\Gamma_e)$ and let the Neumann-Poincar\'e-type operator $\Kbb^* : \Hcal \to \Hcal$ be
defined by
\begin{equation} \label{eq:K*}
 \Kbb^*:= \begin{bmatrix}
 \ds -\Kcal_{\Gamma_i}^* & \ds -\pd{}{\nu_i} \Scal_{\Gamma_e} \\
 \nm \ds \pd{}{\nu_e} \Scal_{\Gamma_i} & \ds \Kcal_{\Gamma_e}^*
 \end{bmatrix},
\end{equation}
and let
\begin{equation} \label{defg}
\Phi:=\begin{bmatrix} \ds \Gvf_i \\ \Gvf_e \end{bmatrix} , \quad
g:= \begin{bmatrix} \ds \frac{\p F}{\p \nu_i} \\ \nm \ds -
\frac{\p F}{\p \nu_e} \end{bmatrix}.
\end{equation}
Then, \eqref{matrixeq3} can be rewritten in the form
\begin{equation}\label{matrixeq4 ep=1}
(  \Ibb^\Gd + \Kbb^* ) \Phi = g,
\end{equation}
where $\Ibb^\Gd$ is given by
\begin{equation}\label{II}
\Ibb^\Gd = \begin{bmatrix}
 \ds z_i^\Gd I  & 0 \\
 0  &  z_e^\Gd I
 \end{bmatrix} .
\end{equation}

%%%%%%%%%%%%%%%%%%%%%%%%%%%%%%%%%%%%
\section{Eigenvalues of the NP operator}
%%%%%%%%%%%%%%%%%%%%%%%%%%%%%%%%%%%%

It is proved in \cite{acklm} that for arbitrary-shaped domains
$\GO$ and $D$ the spectrum of the NP operator $\Kbb^*$ lies in
$[-1/2, 1/2]$, and if $\GO$ and $D$ are concentric disks, the
eigenvalues of $\Kbb^*$ on $\Hcal_0$ are $\pm \rho^n/2$, $n=1,2,
\ldots$. In this section we compute the eigenvalues $\Kbb^*$ on
$\Hcal$ when $\GO$ and $D$ are concentric disks or balls.

%%%%%%%%%%%%%%%%%%%%%%%%%%%%%%%%%%%%
\subsection{Two dimensions}
%%%%%%%%%%%%%%%%%%%%%%%%%%%%%%%%%%%%

Let $\GG=\{ |x|=r_0\}$ in two dimensions. It is known that for each integer
$n$
\begin{equation} \label{Single circular}
\Scal_\Gamma[e^{in\theta}](x) = \begin{cases}
 \ds - \frac{r_0}{2|n|} \left(\frac{r}{r_0}\right)^{|n|}
 e^{in\theta} \quad & \mbox{if } |x|=r < r_0, \\
 \nm
 \ds - \frac{r_0}{2|n|} \left(\frac{r_0}{r}\right)^{|n|}
 e^{in\theta} \quad & \mbox{if } |x|=r > r_0.
 \end{cases}
\end{equation}
Moreover,
 \beq\label{kcal-circ}
 \Kcal_\GG^* [e^{in \theta}] =0 \quad \forall n \neq 0,
 \eeq
and
 \beq\label{kcal-circ-zero}
 \Kcal_\GG[1]= \frac{1}{2}.
 \eeq
In other words, $\Kcal_\GG$ is a rank 1 operator whose only non-zero eigenvalue is $1/2$.

Using \eqnref{kcal-circ}, it is proved that eigenvalues of $\Kbb^*$ on $\Hcal_0$ are $\pm \rho^2/2$ (see \cite{acklm}). We now show that $\pm 1/2$ are also eigenvalues of $\Kbb^*$ on $\Hcal_0$. These eigenvalues are of interest in relation to estimation of stress concentration \cite{ackly}. Using \eqnref{kcal-circ-zero} we have
\begin{equation} \label{Singlezero2D}
\Scal_\GG[1](x) = \begin{cases}
 \ds \log r_0 \quad & \mbox{if } |x|=r < r_0, \\
 \nm
 \ds \log |x| \quad & \mbox{if } |x|=r > r_0,
 \end{cases}
\end{equation}
and hence
 \beq\label{single-cir-nor-zero}
\frac{\p }{\p r} \Scal_\GG[1](x) =  \begin{cases}
 \ds 0 \quad & \mbox{if } |x|=r < r_0, \\
 \nm
 \ds \frac{1}{r} \quad & \mbox{if } |x|=r > r_0.
 \end{cases}
 \eeq
It then follows that
\beq
\Kbb^* \begin{bmatrix} a \\ b \end{bmatrix} =
\begin{bmatrix}
-\frac{1}{2} & 0 \\ \frac{1}{r_e} & \frac{1}{2} \end{bmatrix} \begin{bmatrix} a \\ b \end{bmatrix},
\eeq
where $a$ and $b$ are constants. So $\pm 1/2$ are eigenvalues of $\Kbb^*$.

We summarize our findings in the following proposition.
\begin{prop}
The eigenvalues of $\Kbb^*$ defined on concentric circles in two dimensions are
\beq
-\frac{1}{2}, \ \frac{1}{2}, \ -\frac{1}{2} \rho^{n} , \ \frac{1}{2} \rho^{n}, \ n=1,2,\ldots,
\eeq
and corresponding eigenfunctions are
\beq
\begin{bmatrix} 1 \\ -\frac{1}{r_e} \end{bmatrix}, \ \begin{bmatrix} 0 \\ 1 \end{bmatrix}, \ \begin{bmatrix} e^{\pm in\theta} \\ \rho e^{\pm in\theta} \end{bmatrix}, \ \begin{bmatrix} e^{\pm in\theta} \\ -\rho e^{\pm in\theta} \end{bmatrix}, \ n=1,2,\ldots.
\eeq
\end{prop}

%%%%%%%%%%%%%%%%%%%%%%%%%%%%%%%%%%%%
\subsection{Three dimensions}
%%%%%%%%%%%%%%%%%%%%%%%%%%%%%%%%%%%%

Let $Y^m_n(\hat x)$ $(m=-n,-n+1,\ldots,0, 1,\ldots,n)$ be the
orthonormal spherical harmonics of degree $n$. Here $\hat x=
\frac{x}{|x|}$. Then $|x|^nY^m_n(\hat x)$ is harmonic in $\Bbb
R^3$.
\begin{lem}\label{eigen-sphere}
Let $\GG=\{ |x|=r_0\}$ in three dimensions. We have for
$n=0,1,\ldots$ \beq \label{k_star_ball} \Kcal_\GG^*[Y_n^m](x) =
\frac{1}{2(2n+1)} Y_n^m (\hat x), \quad |x|=r_0, \quad m=-n,
\ldots, n. \eeq
\end{lem}
\pf
It is proved in \cite[Lemma 2.3]{ks99} that
\beq\label{ks99}
\Kcal_\GG^*[\Gvf](x) = -\frac{1}{2r_0} \Scal_\GG[\Gvf](x), \quad |x|=r_0
\eeq
for any function $\Gvf \in L^2(\GG)$. So it follows from \eqnref{singlejump} that
\beq
\frac{\p}{\p r} \Scal_\GG[\Gvf]\big|_{-}(x) + \frac{1}{2r_0} \Scal_\GG[\Gvf](x) = -\frac{1}{2} \Gvf(x), \quad |x|=r_0.
\eeq
Let $\Gvf(x)=Y_n^m(\hat x)$. Then since $\Scal_\GG[Y_n^m](x)$ and $|x|^n Y_n^m(\hat x)$ are harmonic functions in $\{|x| < r_0\}$, we have
\beq\label{scal-harm}
\Scal_{\GG} [Y_n^m](x) = - \frac{1}{2n+1} \frac{r^n}{r_0^{n-1}} Y_n^m (\hat x),\quad \mbox{for}~|x|=r \le r_0,
\eeq
and \eqnref{k_star_ball} follows from \eqnref{ks99}. \qed

Lemma \ref{eigen-sphere} says that the eigenvalues of
$\Kcal_\GG^*$ on $L^2(\GG)$ when $\GG$ is a sphere are
$\frac{1}{2(2n+1)}$, $n=0,1,\ldots$, and their multiplicities are
$2n+1$.

By \eqnref{scal-harm}, we have
\beq
\pd{}{\nu_i} \Scal_{\GG_e}[Y_n^m](x) = - \frac{n}{2n+1} \left( \frac{r_i}{r_e}\right)^{n-1} Y_n^m (\hat x), \quad |x|=r_i.
\eeq
Similarly, we have
$$
\Scal_{\GG_i} [Y_n^m](x) = - \frac{1}{2n+1} \frac{r_i^{n+2}}{r^{n+1}} Y_n^m (\hat x),\quad \mbox{for}~|x|=r \ge r_i,
$$
and hence
\beq
\pd{}{\nu_e} \Scal_{\GG_i}[Y_n^m](x) =\frac{n+1}{2n+1} \left( \frac{r_i}{r_e}\right)^{n+2} Y_n^m (\hat x), \quad |x|=r_e.
\eeq

We now have for constants $a$ and $b$
\beq\label{kstasphere}
\Kbb^*\begin{bmatrix}
a Y_n^m \\ \nm b Y_n^m
\end{bmatrix} = \begin{bmatrix}
\left( - \frac{a}{2(2n+1)} +b\frac{n}{2n+1} \rho^{n-1} \right) Y_n^m \\ \left( a\frac{n+1}{2n+1} \rho^{n+2} +\frac{b}{2(2n+1)}\right) Y_n^m \end{bmatrix} =
\begin{bmatrix}
- \frac{1}{2(2n+1)} & \frac{n}{2n+1} \rho^{n-1} \\ \nm \frac{n+1}{2n+1} \rho^{n+2} & \frac{1}{2(2n+1)} \end{bmatrix} \begin{bmatrix}
a Y_n^m \\ \nm b Y_n^m
\end{bmatrix}.
\eeq Thus we have the following result.
\begin{prop}\label{sphere}
The eigenvalues of $\Kbb^*$ defined on two concentric spheres are
\beq
\pm \frac{1}{2(2n+1)}\sqrt{1+4n(n+1)\rho^{2n+1}}, \quad n=0,1,\ldots,
\eeq
and corresponding eigenfunctions are
\begin{align}
\begin{bmatrix} (\sqrt{1+4n(n+1)\rho^{2n+1}}-1)Y_n^m\\ 2(n+1)\rho^{n+2}Y_n^m \end{bmatrix}, \quad \begin{bmatrix} (-\sqrt{1+4n(n+1)\rho^{2n+1}}-1)Y_n^m\\ 2(n+1)\rho^{n+2}Y_n^m \end{bmatrix}, \ m= -n,\ldots,n,
\end{align}
respectively.
\end{prop}

It is quite interesting to observe that if we let $\frac{1}{2}= \Gl_0 \ge \Gl_1 \ge \ldots$ be the eigenvalues of $\Kcal_\GG$ for a disk or a sphere enumerated according to their multiplicities, then the eigenvalues $\mu_n$ of $\Kbb^*$ satisfy
\beq
\mu_n = \pm \Gl_n + O(\rho^n).
\eeq

%%%%%%%%%%%%%%%%%%%%%%%%%%%%%%%%%%%%
\section{Anomalous localized resonance in two dimensions}
%%%%%%%%%%%%%%%%%%%%%%%%%%%%%%%%%%%%

In this section we consider the CALR when the domains $\GO$ and
$D$ are concentric disks. We first observe that $z_i^\Gd$ and
$z_e^\Gd$ converges to non-zero numbers as $\Gd$ tends to $0$ if
$\Ge_c \neq- \Ge_s \neq 1$. So, in this case CALR does not occur
regardless of the location of the source. Furthermore, if
$\Ge_c=\Ge_s=1$, a thorough study was done in \cite{acklm}. It is
proved in \cite{acklm} that if the source $f$ is supported inside
the critical radius $r_*=\sqrt{r_e^3/r_i}$, then the weak CALR
occurs, namely, \beq \limsup_{\Gd \to 0} E_\Gd = \infty. \eeq
Moreover, if $F$ is the Newtonian potential of $f$ and the Fourier
coefficients $g_e^n$ of $-\pd{F}{\nu_e}$ satisfies the following
gap property:
\begin{quote}
[GP] There exists  a sequence $ \{n_k\} $ with $|n_1|<|n_2|<\cdots$ such that
$$
\lim_{k \to \infty} \rho^{|n_{k+1}|-|n_k|}\frac{|g_e^{n_k}|^2}{|n_k|\rho^{|n_k|}} = \infty,
$$
\end{quote}
then CALR occurs, namely \beq \lim_{\Gd \to 0} E_\Gd = \infty,
\eeq
and $V_\Gd/\sqrt{E_\Gd}$ goes to zero outside the radius $\sqrt{r_e^3/r_i}$.

The remaining two cases are when $\Ge_c \neq -\Ge_s=1$ and $\Ge_c = -\Ge_s \neq 1$. In these cases, we have the following theorem.
\begin{thm}\label{thm41}
\begin{itemize}
\item[(i)] If $\Ge_c = -\Ge_s \neq 1$, then CALR does not occur,
{\it i.e.},
\beq
E_\Gd \le C
\eeq
for some $C>0$. (We note, however, that there will be CALR for appropriately placed sources inside the core, as
can be seen from the fact that the equations are invariant under conformal transformations, and in particular under the inverse
transformation $1/z$ where $z=x_1+i x_2$, which in effect interchanges the roles of the matrix and core.)

\item[(ii)] If $\Ge_c \neq- \Ge_s = 1$, then weak CALR occurs and
the critical radius is $r_*= r_e^2 r_i^{-1}$, {\it i.e.}, if the
source function is supported inside $r_*$ (and its Newtonian
potential does not extend harmonically to $\Rbb^2$), then \beq
\limsup_{\Gd \to 0} E_\Gd = \infty, \eeq
and there exists a constant $C$ such that
\beq\label{bounded1}
|V_\delta(x)|<C
\eeq
for all $ x$ with $|x|>r_e^3/r_i^2$.
\item[(iii)] In addition to the assumptions of (ii), the Fourier coefficients $g_e^n$ of $-\pd{F}{\nu_e}$ satisfies the following gap property:
\begin{quote}
[GP2] There exists  a sequence $ \{n_k\} $ with $|n_1|<|n_2|<\cdots$ such that
$$
\lim_{k \to \infty} \rho^{2(|n_{k+1}|-|n_k|)}\frac{|g_e^{n_k}|^2}{|n_k|\rho^{|n_k|}} = \infty,
$$
\end{quote}
then the CALR occurs, {\it i.e.}, \beq \lim_{\Gd \to 0} E_\Gd =
\infty, \eeq
and $V_\Gd/\sqrt{E_\Gd}$ goes to zero outside the radius $r_e^3/r_i^2$, as implied by \eqnref{bounded1}.
\end{itemize}
\end{thm}

Before proving Theorem \ref{thm41} we make a remark on the Gap Properties [GP] and [GP2]. One can easily see that [GP] is weaker than [GP2], namely, if [GP] holds, so does [GP2].

The rest of this section is devoted to the proof of Theorem \ref{thm41}. As was proved in \cite{acklm}, we have
\begin{align*}
\frac{\p}{\p \nu_e} \Scal_{\GG_i}[e^{in\theta}](x) &=  \frac{1}{2} \rho^{|n|+1}
 e^{in\theta}, \\
\frac{\p}{\p \nu_i} \Scal_{\GG_e}[e^{in\theta}](x) &=  -\frac{1}{2} \rho^{|n|-1}
 e^{in\theta}.
\end{align*}
Using these identities, one can see that
if $g$ defined by \eqnref{defg} has the Fourier series expansion
 $$
 g= \sum_{n \neq 0} \begin{bmatrix} g_i^n\\ g_e^n \end{bmatrix} e^{in\theta},
 $$
then the integral equations \eqnref{matrixeq4 ep=1} are equivalent to
\beq\label{inteqn2}
\begin{cases}
\ds z_i^\delta \Gvf_i^n + \frac{\rho^{|n|-1}}{2} \Gvf_e^n =  g_i^n , \\
\nm
\ds z_e^\delta \Gvf_e^n + \frac{\rho^{|n|+1}}{2} \Gvf_i^n =  g_e^n
\end{cases}
\eeq
for every $|n|\geq 1$. It is readily seen that the solution
$\Phi=(\Gvf_i,\Gvf_e)$ to \eqnref{inteqn2} is given by
\begin{align*}
\Gvf_i &= 2 \sum_{n\neq 0} \frac{2z_e^\delta g_i^n-\rho^{|n|-1}g_e^n} {4z_i^\delta z_e^\delta-\rho^{2|n|}}  e^{i n
\theta}, \\
\Gvf_e &= 2 \sum_{n\neq 0} \frac{2z_i^\delta g_e^n-\rho^{|n|+1}g_i^n} {4z_i^\delta z_e^\delta-\rho^{2|n|}}  e^{i n
\theta}.
\end{align*}

If the source is located outside the structure, {\it
i.e.}, $f$ is supported in $|x| > r_e$, then
the Newtonian potential of $f$, $F$, is harmonic in $|x| \le r_e$ and
\beq\label{Ffourier}
F(x)=c-\sum_{n\neq 0}
\frac{g_e^n}{|n|r_e^{|n|-1}}r^{|n|}e^{in\theta}, \quad |x| \le
r_e.
\eeq
Thus we have
 \beq\label{gin}
 g_i^n = - g_e^n \rho^{|n|-1}.
 \eeq
So we have
\begin{equation}\label{formulaofphi}
\begin{cases}
\ds\Gvf_i = -2 \sum_{n\neq 0} \frac{(2z_e^\delta +1)\rho^{|n|-1}g_e^n} {4z_i^\delta z_e^\delta-\rho^{2|n|}}  e^{i n
\theta},\\
\ds\Gvf_e = 2 \sum_{n\neq 0} \frac{(2z_i^\delta +\rho^{2|n|})g_e^n} {4z_i^\delta z_e^\delta-\rho^{2|n|}}  e^{i n
\theta}.
\end{cases}
\end{equation}
Therefore, from \eqref{Single circular} we find that
\beq\label{SSoutside} \Scal_{\Gamma_i}[\Gvf_i](x) +
\Scal_{\Gamma_e}[\Gvf_e](x) = \sum_{ n\neq 0} \frac{2(r_i^{2|n|}z_e^\delta -
r_e^{2|n|}z_i^\delta)}{|n|r_e^{|n|-1}(4z_i^\delta z_e^\delta - \rho^{2|n|})}
\frac{g_e^n}{r^{|n|}} e^{i n \theta} , \quad r_e < r=|x|, \eeq and
\begin{align}
\Scal_{\Gamma_i}[\Gvf_i](x)
& =  -\sum_{ n\neq 0} \frac{ r_i^{2|n|}(2z_e^\delta+1)}{|n|r_e^{|n|-1}(\rho^{2|n|}-4z_i^\delta z_e^\delta)}
\frac{g_e^n}{r^{|n|}} e^{i n \theta} ,  \quad r_i < r=|x| < r_e, \label{singleexp1} \\
\Scal_{\Gamma_e}[\Gvf_e](x) & = \sum_{ n\neq 0}
\frac{(2z_i^\delta+\rho^{2|n|})}{|n|r_e^{|n|-1}(\rho^{2|n|}-4z_i^\delta z_e^\delta)}
g_e^n r^{|n|} e^{i n \theta}, \quad r_i < r=|x| < r_e.
\label{singleexp2}
\end{align}

We obtain the following lemma.

\begin{lem} \label{lemest}
There exists $\delta_0$ such that
\beq\label{essest}
E_\Gd \approx \begin{cases} \ds\sum_{n \ne 0}
\frac{\delta|g_e^n|^2}{|n|(\delta^2+\rho^{4|n|})},\quad \mbox{if }  \Ge_c \neq \Ge_s=1 ,\\
\ds\sum_{n \ne 0}
\frac{\delta\rho^{2|n| }|g_e^n|^2}{|n|(\delta^2+\rho^{4|n|})},\quad \mbox{if }  \Ge_c = \Ge_s  \ne 1, \end{cases} \eeq uniformly
in $\delta \le \delta_0$.
\end{lem}

\pf \, Using \eqnref{Ffourier}, \eqnref{singleexp1}, and
\eqnref{singleexp2}, one can see that
$$
V_\delta (x) =  c+ r_e \sum_{n\neq 0} \left[
\frac{r_i^{2|n|}}{r^{|n|}}  -2z_i^\delta r^{|n|}
\right] \frac{(2z_e^\delta+1)g_e^n e^{i n \theta}}{|n|r_e^{|n|}(4z_i^\delta z_e^\delta -
\rho^{2|n|})} .
$$
We check that
\begin{align*}
\left| \nabla\left( \left(\frac{r_i^{2|n|}}{r^{|n|}} - 2 z_i^\delta r^{|n|}\right) e^{in\theta} \right)\right|^2 = \frac{2|n|^2}{r^2} \left| \frac{r_i^{2|n|}}{r^{|n|}} - 2 z_i^\delta r^{|n|}\right|^2.
\end{align*}
Then straightforward computations yield that
\begin{align*}
\int_{B_e\setminus B_i} \delta |\nabla V_\delta|^2& \approx
\sum_{n\ne0} \delta \left|\frac{2z_e^\delta +1}
{4z_i^\delta z_e^\delta-\rho^{2|n|}}\right|^2 (4|z_i^\delta|^2 + \rho^{2|n|})
\frac{|g_e^n|^2}{|n|} .
\end{align*}
If $\delta$ is sufficiently small, then one can also easily show
that
 $$
 |4z_i^\delta z_e^\delta-\rho^{2|n|}| \approx \delta + \rho^{2|n|}.
 $$
Therefore we get \eqnref{essest} and the proof is complete.
\qed

First, if $\Ge_c=-\Ge_s \ne 1$, then
 $$
 E_\Gd \approx \sum_{n \ne 0}
\frac{\delta\rho^{2|n| }|g_e^n|^2}{|n|(\delta^2+\rho^{4|n|})} \le \sum_{n \ne 0}
\frac{|g_e^n|^2}{2|n|} \le \frac{1}{2} \left\| \pd{F}{\nu_e} \right\|_{L^2(\GG_e)} \le C \| f \|_{L^2(\Rbb^2)}.
$$

Suppose that $\Ge_c \ne-\Ge_s = 1$, and let
 \beq
 N_\delta = \frac{\log \delta}{2\log \rho}.
 \eeq
If $|n| \le N_\delta$, then $\delta \le \rho^{2|n|}$, and hence
 \beq \label{easier009}
\sum_{n \ne 0} \frac{\delta|g_e^n|^2}{|n|(\delta^2+\rho^{4|n|})} \ge  \sum_{0 \ne |n| \le N_\delta} \frac{\delta|g_e^n|^2}{|n|(\delta^2+\rho^{4|n|})} \ge \sum_{0 \ne |n| \le N_\delta} \frac{\delta|g_e^n|^2}{|n|\rho^{4|n|}}.
 \eeq
If the following holds\beq\label{limsup-cond}
\limsup_{n\rightarrow\infty} \frac{|g_e^n|^2}{|n|\rho^{2|n|}}=
\infty, \eeq then one can show as in \cite{acklm} that there is a
sequence $\{|n_k|\}$ such that \beq \lim_{k \to \infty}
E_{\rho^{|n_k|}} = \infty. \eeq

Suppose that the source function $f$ is supported inside the
critical radius $r_*= r_e^2r_i^{-1}$ (and outside $r_e$). Then its
Newtonian potential $F$ cannot be extended harmonically in $|x| <
r_*$ in general. So, if $F$ is given by
 \beq\label{fser}
 F = c- \sum_{n\ne 0} a_n r^{|n|} e^{in\theta},\quad r< r_e + \Ge
 \eeq
for some $\Ge >0$, then the radius of convergence of the series is
less than $r_*$. Thus we have
 \beq\label{easier22}
 \limsup_{|n|\rightarrow\infty} |a_n|^2 r_*^{2|n|}= \infty.
 \eeq
Since $g^n_e= |n| a_n r_e^{|n|-1}$, \eqnref{limsup-cond} holds.

By \eqnref{SSoutside}, we know
\beq
|V_\delta| \le |F| + C\sum_{n\ne 0} \frac{r_e^{|n|}}{\delta+\rho^{2|n|}}\frac{|g_e^n|}{r^{|n|}}
\le |F| + C\sum_{n\ne 0} \frac{r_e^{3|n|}}{r_i^{2|n|}}\frac{|g_e^n|}{r^{|n|}}\le C'  \eeq
if $r> r_e^3/r_i^{2}$.
Thus (ii) is proved.

We now prove (iii). We emphasize that [GP2] implies \eqnref{limsup-cond}, but the converse may not be true. On the other hand [GP2] holds if
\beq\label{lim-cond}
\lim_{n\rightarrow\infty} \frac{|g_e^n|^2}{|n|\rho^{2|n|}}= \infty.
\eeq
So we may regard the condition [GP2] something between \eqnref{limsup-cond} and \eqnref{lim-cond}.

Suppose that [GP2] holds. If we take $\delta=\rho^{2\alpha}$ and let $k(\alpha)$ be the number such that $$ |n_{k(\alpha)}| \le \alpha < |n_{k(\alpha)+1}|,$$  then
 \beq
 \sum_{0 \ne |n| \le N_{\delta}} \frac{\delta|g_e^n|^2}{|n|\rho^{4|n|}} = \rho^{2\alpha} \sum_{0 \ne |n| \le \alpha} \frac{|g_e^n|^2}{|n|\rho^{4|n|}}\ge \rho^{2(|n_{k(\alpha)+1}|-|n_{k(\alpha)}|)}\frac{|g_e^{n_{k(\alpha)}}|^2}{|n_{k(\alpha)}|\rho^{2|n_{k(\alpha)}|}} \rightarrow \infty,
 \eeq
 as $\alpha\to\infty$. Combined with Lemma \ref{lemest} and \eqnref{easier009}, it gives us (iii).

%%%%%%%%%%%%%%%%%%%%%%%%%%%%%%%%%%%%%%%%
\section{Non-occurrence of CALR in 3D}
%%%%%%%%%%%%%%%%%%%%%%%%%%%%%%%%%%%%%%%%

In this section we show that CALR does not occur in a radially symmetric three
dimensional coated sphere structure when the core, matrix and shell are isotropic. We have the following theorem.
\begin{thm}
Suppose that $\GG_e$ and $\GG_i$ are concentric spheres. For any
$\Ge_c$ and $\Ge_s$, there is a constant $C$ independent of $\Gd$
such that if $V_\delta$ is the solution to \eqnref{basiceqn}, then
\beq  \label{main3D} \int_{\GO \setminus D} \Gd |\nabla V_\Gd|^2
\le C \| f \|_{L^2(\Rbb^3)}^2. \eeq
\end{thm}
\pf
Suppose that $\pd{}{\nu_e} F$ has the Fourier series expansion
\beq
\pd{}{\nu_e} F = -\ds\sum_{n=0}^\infty \sum_{m=-n}^n g_{mn}^e Y_m^n.
\eeq
Then one can show as in \eqnref{gin} that
\beq
\pd{}{\nu_i} F = -\ds\sum_{n=0}^\infty \sum_{m=-n}^n g_{mn}^e \rho^{n-1} Y_m^n.
\eeq
By solving the integral equation \eqnref{matrixeq3} using \eqnref{kstasphere}, we obtain
\begin{align}
\Gvf_i &= -\sum_{n=0}^\infty \sum_{m=-n}^n \GD_n^{-1}\rho^{n-1} \left( z_e^\Gd + \frac{1}{2} \right) g_{mn}^e Y_n^m,\\
\Gvf_e &= \sum_{n=0}^\infty \sum_{m=-n}^n \GD_n^{-1}\rho^{n-1} \left( z_i^\Gd -\frac{1}{2(2n+1)}+ \frac{n+1}{2(2n+1)}\rho^{2n+1} \right) g_{mn}^e Y_n^m,
\end{align}
where
$$
\GD_n:= \left(z_i^\Gd - \frac{1}{2(2n+1)}\right)\left(z_e^\Gd + \frac{1}{2(2n+1)}\right)+ \frac{n(n+1)}{(2n+1)^2} \rho^{2n+1}.
$$

Suppose for simplicity that $\Ge_c=-\Ge_s=1$, so that $z_i^\Gd$ and
$z_e^\Gd$ given by \eqnref{lam e mu} simplify to
$$
z_i^\Gd = z_e^\Gd = \frac{i\Gd}{2(2-i\Gd)}.
$$
Then one can see that if $\Gd$ is sufficiently small, then
$$
|\GD_n| \approx \Gd^2 + n^{-2}.
$$
So we have
\begin{align*}
\Gd \| \Gvf_i \|_{L^2(\GG_i)}^2 &\le C \sum_{n=0}^\infty \sum_{m=-n}^n \frac{\Gd \rho^{2n}}{(\Gd^2 + n^{-2})^2} |g_{mn}^e|^2 \\
& \le C \sum_{n=0}^\infty \sum_{m=-n}^n n^3 \rho^{2n} |g_{mn}^e|^2 \\
& \le C \sum_{n=0}^\infty |g_{mn}^e|^2 \le C \| f \|_{L^2(\Rbb^3)}^2.
\end{align*}
and
$$
\Gd \| \Gvf_e \|_{L^2(\GG_e)}^2 \le C \sum_{n=0}^\infty \sum_{m=-n}^n \frac{\Gd \rho^{2n}}{\Gd^2 + n^{-2}} |g_{mn}^e|^2 \le C \| f \|_{L^2(\Rbb^3)}^2.
$$
Therefore we have
\begin{align*}
\int_{\GO \setminus D} \Gd |\nabla V_\Gd|^2 &= \int_{\GO \setminus D} \Gd |\nabla F|^2 + \int_{\GO \setminus D} \Gd |\nabla ( \Scal_{\GG_i}[\Gvf_i] + \Scal_{\GG_e}[\Gvf_e])|^2 \\
&\le \int_{\GO \setminus D} \Gd |\nabla F|^2 + \int_{\GO \setminus D} \Gd |\nabla ( \Scal_{\GG_i}[\Gvf_i] + \Scal_{\GG_e}[\Gvf_e])|^2 \\
&\le \int_{\GO \setminus D} \Gd |\nabla F|^2 + \Gd (\| \Gvf_i \|_{L^2(\GG_i)}^2 + \| \Gvf_e \|_{L^2(\GG_e)}^2) \le C \| f \|_{L^2(\Rbb^3)}^2.
\end{align*}
If $\Ge_s\neq -1$ or/and $\Ge_c\neq 1$, then the same argument can
be applied to arrive at \eqnref{main3D}. This completes the proof.
\qed

\end{document}